\documentclass[11pt, fleqn, leqno]{article}
\usepackage{graphicx}
\usepackage{amsmath}
\usepackage{amsfonts}
\usepackage{amsbsy}
\usepackage{latexsym}
\numberwithin{equation}{section}

\newcommand{\be}{\begin{equation}}
\newcommand{\ee}{\end{equation}}
\newcommand{\ba}{\begin{array}}
\newcommand{\ea}{\end{array}}

\newcommand{\dv}{\frac{\v^{\prime}(z)-\v^{\prime}(y)}{z-y}}

\renewcommand{\v}{{\mathsf{v}}}
\newcommand{\M}{{\mathsf{M}}}
\newcommand{\U}{{\mathsf{U}}}

\newcommand{\pr}{\prime}
\newcommand{\bt}{\beta}

\renewcommand{\b}{\beta}
\renewcommand{\l}{\lambda}

\newcommand{\bea}{\begin{eqnarray}}
\newcommand{\eea}{\end{eqnarray}}

\newcommand{\al}{\alpha}

\begin{document}
\newtheorem{pro}[thm]{Proposition}
\newtheorem{lem}[thm]{Lemma}
\newtheorem{cor}[thm]{Corollary}
\newcommand{\C}{{\mathsf{C}}}
\newcommand{\A}{{\mathsf{A}}}
\newcommand{\G}{{\mathsf{\Gamma}}}
\renewcommand{\a}{{\mathsf{a}}}

\title{Jacobi Polynomials from Compatibility Conditions
\thanks{Research supported by NSF
grant DMS 99-70865 and by EPSRC grant GR/S14108}}

\author{Yang Chen\\
Department of Mathematics \\
Imperial College\\
180 Queen's Gate, \\London  SW7 2BZ, UK\\
email:  ychen@imperial.ac.uk  \\
\and Mourad Ismail  \\ Department of Mathematics\\
University of South Florida, \\Tampa,  
Florida 33620 \\USA\\
\; \;  email: ismail@math.usf.edu }
\date{}
\maketitle

\begin{abstract} 
\noindent
We revisit the ladder operators for orthogonal polynomials
and re-interpret two supplementary conditions as compatibility
conditions of two linear over-determined systems; 
one involves the variation of the polynomials with respect to 
the variable $z$ (spectral parameter) and the other a 
recurrence relation in $n$ (the lattice variable). For the 
Jacobi weight $w(x)=(1-x)^{\al}(1+x)^{\bt},\;x\in[-1,1],$ we 
show how to  use the compatibility conditions to explictly
determine  the recurrence coefficients of the monic Jacobi
polynomials.  
 
\end{abstract}

\bigskip
\noindent
{\bf Running title}: Jacobi polynomials

%\newpage 

\bigskip

\section{Introduction and Preliminaries.}

We begin with some notation. 
Let $P_n(x)$ be monic polynomials of degree $n$ in $x$ 
and orthogonal, with respect to a weight, $w(x),\;x\in[a,b];$ 
\bea
\int_{a}^{b}P_m(x)P_n(x)w(x)dx=h_n\delta_{m,n}.
\eea
% For convenience, we suppose $w(a)=w(b)=0,$ 
% (although this is not necessary), 
We further assume that  
$\v^{\prime}(z):=-w^{\prime}(z)/w(z)$ exists and that\\
\noindent
$y^n [\v'(x)-\v'(y)] w(y)/(x-y)$ is integrable 
on $[a,b]$ for
all $n$, 
$n =0, 1. \cdots$.  From the orthogonality condition there
follows the recurrence  relation,
\bea
zP_n(z)=P_{n+1}(z)+\al_nP_n(z)+\bt_nP_{n-1}(z),
n=0,1,...,\;
\eea
where $\bt_0 P_{-1}(z):=0$, $\al_n,\;n=0,1,2,...$ is real and
$\bt_n>0, \;n=1,2,...$

In this paper we
describe a formalism  which derives  
properties of  orthogonal 
polynomials, and  their recurrence coefficients,
from the knowledge of the weight function. We believe this is
a new and interesting approach to orthogonal
polynomials. In order to keep this work accessible we will only
include the example of Jacobi polynomials. We defer in a future 
publication, the analysis in
the case of the generalized Jacobi weights \cite{Ne},  
\cite{Is2}. In the Jacobi case we find the recurrence relations
in \S 2. In \S 3 we show how our approach leads to the
evaluation of monic Jacobi polynomials at $x = \pm 1$.  We
also show that the evaluation of a Jacobi polynomial at $x=1$ 
or $x = -1$ leads to explicit representations of the Jacobi
polynomials. Closed form expressions for the normalization
constants $h_n$ are also found.  
 
The actions of the ladder operators on $P_n(z)$ and
$P_{n-1}(z)$  are,
\bea
\mbox{} \left(\frac{d}{dz}+B_n(z)\right)P_n(z) =
\bt_nA_n(z)P_{n-1}(z)
\\
\mbox{} \left(\frac{d}{dz}-B_n(z)-\v^{\pr}(z)\right)P_{n-1}(z)
= -A_{n-1}(z)P_n(z),
\eea
with 
\bea
\mbox{}A_n(z)&:=& \left.\frac{w(y)\, P_n^2(y)}
{h_n(y -z)}\right|_{y=a}^{y=b}
+ \frac{1}{h_n} \int_{a}^{b} \dv P_n^2(y)w(y)dy,
\\
\mbox{}B_n(z)&:=& \left.\frac{w(y)\, P_n(y)P_{n-1}(y)}  
{h_n(y -z)}\right|_{y=a}^{y=b} \\
&{}& \qquad + \frac{1}{h_{n-1}} \int_{a}^{b}  \dv
P_{n-1}(y)P_n(y) w(y)dy, \nonumber 
% \end{gathered}
\eea
where we have used the supplementary condition,
$$
B_{n+1}(z)+B_n(z)=(z-\al_n)A_n(z)-\v^{\prime}(z),\eqno(S_1)
$$
to arrive at (1.4). The equations (1.3)--(1.6) and the 
supplementary condition $(S_1),$
was derived by Bonan and Clark \cite{Bo:Cl}, 
 Bauldry \cite{Ba}, and Mhaskar \cite{Mh}  for polynomial $\v,$ and
the Authors
\cite{Ch:Is}  for  general $\v.$ Ismail and Wimp \cite{Is:Wi}
identified the additional supplementary condition,
$$
B_{n+1}(z)-B_{n}(z)=\frac{\bt_{n+1}A_{n+1}(z)-\bt_nA_{n-1}(z)-1}
{z-\al_n}. \eqno(S_2)
$$

Our thesis in this work is that the supplementary conditions,
$(S_1)$ and
$(S_2),$ being identities in $n, (n >0)$ and $z\in{\bf C}\cup\infty$
 have the
information needed to determine the recurrence
coefficients and other auxilliary quantities. We illustrate
this by systematically using $(S_1)$ and $(S_2)$ to
determine most of the properties of the Jacobi polynomials. See 
\cite{Sz}, \cite{An:As:Ro}, and \cite{Rai}, for information concerning the
Jacobi polynomials.  In describing our results  
we shall follow the standard notation for shifted factorial
and hypergeometric functions in \cite{An:As:Ro}, \cite{Rai}.

Below, we reinterpret $(S_2).$
We set 
\bea
\Phi_n(z)&:=&\left(\begin{matrix}P_n(z) \\
                         P_{n-1}(z) \end{matrix}\right),
\\
\M_n(z)&:=&\left(\begin{matrix} -B_n(z)&\bt_nA_n(z)
\\                      
-A_{n-1}(z)&B_n(z)+\v^{\pr}(z)\end{matrix}\right),
\\
\U_n(z)&:=&\left(\begin{matrix}z-\al_n &-\bt_n \\
                      1&0 \end{matrix}\right).
\eea
Now equations (1.3) and (1.4) become,
\bea
\Phi^{\pr}_n(z)=\M_n(z)\Phi_n(z),
\eea
and the recurrence relations become,
\bea
\Phi_{n+1}(z)=\U_n(z)\Phi_n(z).
\eea
We find, by requiring (1.10) and (1.11) be compatible; 
\bea
\mbox{} \Phi^{\pr}_{n+1}(z)&=&\M_{n+1}(z)\Phi_{n+1}(z)
\nonumber\\
\mbox{} &=&\M_{n+1}(z)\U_n(z)\Phi_n(z)\nonumber.
\eea
On the other hand 
\bea
\mbox{}
\mbox{} \Phi^{\pr}_{n+1}(z)&=&\U^{\pr}_n(z)\Phi_n(z)+
\U_n(z)\Phi^{\pr}_{n}(z)\nonumber\\
&=&\U^{\pr}_n(z)\Phi_n(z)+\U_n(z)\M_n(z)\Phi_n(z).\nonumber
\eea
We now write the above equations in matrix form as
\bea
\mathsf{S}_n(z)\Phi_n(z)=0,
\eea
where $\mathsf{S}_n(z)$ is the matrix whose entries are 
\bea
\mbox{} \mathsf{S}_n(z)&:=&
\U_n^{\pr}(z)+\U_n(z)\M_n(z)-\M_{n+1}(z)\U_n(z)
\nonumber\\
\mbox{} \mathsf{S}_n^{11}(z) &=&1+(z-\al_n)(B_{n+1}(z)-B_n(z)) \\
&{}& \qquad +
\bt_nA_{n-1}(z)-\bt_{n+1}A_{n+1}(z)
\nonumber\\
\mbox{} \mathsf{S}_n^{12}(z)&=&-\bt_n\left(B_{n+1}(z)
+B_n(z)+\v^{\pr}(z)-(z-\al_n)A_n(z)\right)
\nonumber\\
\mbox{} \mathsf{S}_n^{21}(z)&=&S_n^{12}(z)/\bt_n
\nonumber\\
\mbox{} \mathsf{S}_n^{22}(z)& = & 0. \nonumber
\eea
Here $n=1,2,...$ and $z\in{\bf C}\cup\infty.$
Observe that with $(S_1),$ $S^{12}_n(z) = S^{21}_n(z)=0$. 
This leaves
$S^{11}_n(z) P_n(z)=0$. Since $P_n(z)$ 
does not vanish identically, we must have $S^{11}_n(z)=0,$
which is $(S_2).$ 
It is clear from (1.5) and (1.6) that,
if $\v^{\prime}(z)$ is a rational function then
$A_n(z)$ and
$B_n(z)$ are also rational functions. This is particularly
useful for our purpose,  which is to determine the recurrence
coefficients,
$\al_n$  and $\bt_n.$ In the next section, we illustrate the
method by considering the  Jacobi weight $w^{(\al, \b)}(x) =
(1-x)^\al(1+x)^\b$ for $x \in [-1,1]$. 

Recall that  the numerator polynomials \cite{Sh:Ta}, \cite{Ak}
are 
\bea
Q_n(z):= \int_{-\infty}^\infty\frac{P_n(z)-P_n(y)}{z-y}\;
w(y)\; dy,  
\eea
and $\{P_n(z)\}$ and $\{Q_n(z)\}$ form a basis of solutions of
the recurrence relation. We shall also use the notation 
\bea
F(z) = \int_{-\infty}^\infty\frac{w(y)}{z-y}\; dy, 
\eea
for the Stieltjes transform of the weight function. 

\setcounter{thm}{0}

\setcounter{equation}{0}

\section{Jacobi Weight}
 \setcounter{section}{2}

The Jacobi weight is
$w^{(\al, \b)}(x)=(1-x)^{\al}(1+x)^{\bt};x\in[-1,1],$  
and for now we take
$\al$ and $\bt$ to be strictly positive. It will become clear, using
a real analyticity argument, the  results that follows are also valid
for $\al,\bt > -1.$ Let
$\{{\cal P}_n^{(\al, \b)}(x)\}$ and $\{{\cal Q}_n^{(\al,
\b)}(x)\}$ denote the monic Jacobi polynomials, and their
numerators, respectively, see (1.14). Moreover in the
present example, the Stieltjes transform of $w^{\al,\bt}$
will be denoted by $F^{(\al, \b)}(z)$. 

From (1.5)--(1.6) we find 
\bea
h_n A_n(z) &=& \frac{\al}{1-z}\int_{-1}^1  [{\cal P}_n^{(\al.
\b)}(y)]^2  (1-y)^{\al-1}(1+y)^\b \; dy \nonumber \\
&{}&  +  \frac{\b}{1+z}\int_{-1}^1
 [{\cal P}_n^{(\al. \b)}(y)]^2 (1-y)^{\al}(1+y)^{\b-1} 
\; dy. \nonumber 
\eea
Through  integration by parts, it readily follows that,
\bea
\mbox{} A_n(z)= -\frac{R_n}{z-1}+\frac{R_n}{z+1},
\eea
for some constant $R_n$. Similarly we find 
\bea
 B_n(z)=-\frac{n+r_n}{z-1}+\frac{r_n}{z+1}.
\eea
Here $R_n$ and $r_n$ are given by
\bea
\mbox{} R_n=R_n(\al,\bt)&:=&\frac{\bt}{h_n}
\int_{-1}^{1}\frac{[{\cal P}_n^{(\al,\b)}(y)]^2}{1+y}\;
w^{(\al,\b)}(y)\; dy,\\
\mbox{} r_n = r_n(\al,\bt)&:=&\frac{\bt}{h_{n-1}}
\int_{-1}^{1}\frac{{\cal P}_n^{(\al,\b)}(y) {\cal P}_{n-1}^
{(\al,\b)}(y)}
{1+y}\; w^{(\al,\b)}(y) \, dy.
\eea
  
It is easy to see that
\bea
\mbox{}
R_n(\al,\b)&=& \frac{\bt}{h_n}{\cal P}_n^{(\al,
\b)}(-1)\left[{\cal Q}_n^{(\al, \b)}(-1) 
- F^{(\al, \b)}(-1) P_n^{(\al, \b)}(-1)
\right], \nonumber\\
\mbox{}
r_n (\al, \b) &=& \frac{\bt}{h_{n-1}}{\cal P}_{n-1}^{(\al,
\b)}(-1)
\left[{\cal Q}^{(\al, 
\b)}_n(-1)-F^{(\al, \b)}(-1) {\cal P}^{(\al, \b)}_n
(-1)\right]. \nonumber
\eea
The reader may ask,
``What is the point of this formalism? Since in the attempt to find
$\al_n$ and $\bt_n,$ two new unknown quantities, $R_n$ and $r_n,$
have been  introduced.'' However when $\v'$ is a rational function 
both sides of $(S_1)$ and
$(S_2),$ are rational functions and  by equating coefficients and
residues of  both sides of $(S_1)$ and $(S_2),$
we  shall arrive at four equations which should be sufficient for
the  determination of $R_n$ and $r_n$ as well as $\al_n$ and 
$\bt_n.$ Equating residues at $z=-1$ and $z=+1,$ of $(S_1),$ gives
\bea
\mbox{} -2n-1-r_n-r_{n+1}&=&\al-R_n(1-\al_n)\\
\mbox{} r_n+r_{n+1}&=&\bt-R_n(1+\al_n).
\eea
Similarly, from $(S_2),$ we obtain 
\bea
\mbox{} (r_n-r_{n+1}-1)(1-\al_n)&=&\bt_nR_{n-1}-\bt_{n+1}R_{n+1}
\\
\mbox{} (r_n-r_{n+1})(1+\al_n)&=&\bt_{n+1}R_{n+1}-\bt_nR_{n-1}.
\eea
Observe that $R_n$ can be obtained immediately by adding (2.5) 
and (2.6);
\bea
R_n=\frac{1}{2}(\al+\bt+2n+1).
\eea
The sum (2.7) of (2.8) gives,
\bea
1-\al_n=2(r_n-r_{n+1}),
\eea
while (2.8) minus (2.7) and with (2.9) gives,
\bea
\bt-\al-2n-1-(\al+\bt+2n+1)\al_n=2(r_n+r_{n+1}).
\eea
Now, (2.10) plus (2.11) implies,
\bea
4r_n=\bt-\al-2n-(\al+\bt+2n+2)\al_n
\eea
and (2.11) minus (2.10) implies,
\bea
4r_{n+1}=\bt-\al-2n-2-(\al+\bt+2n)\al_n.
\eea
When (2.12) and (2.13) are made compatible, we obtain a first order 
difference equation satisfied by
$\al_n:$
\bea
\al_{n+1}(\al+\bt+2n+4)-\al_n(\al+\bt+2n)=0,
\eea
which has a very simple ``integrating factor,'' 
$\al+\bt+2n+2$.  Using this, we find,
\bea
\mbox{} \al_n = \frac{C_1}{(2R_n-1)(2R_n+1)}. 
\nonumber
\eea
where $C_1$ is a ``integration'' constant, 
determined by the initial condition,
\bea
\mbox{}\al_0 = \frac{\mu_1}{\mu_0}=\frac{\bt-\al}{\al+\bt+2},\qquad 
C_1 = \bt^2-\al^2.\nonumber
\eea
Here $\mu_j:=\int_{-1}^{1}t^jw(t)dt,\;j=0,1,\cdots$
are the moments.  Therefore we have established 
\bea
\al_n = \frac{\b^2 -\al^2}{(\al+\bt+2n)(\al+\bt+2n+2)}.
\eea

Going back to (2.8) and using (2.10), we see 
that $\bt_n$ satisfies the linear difference equation:
\bea
\bt_{n+1}R_{n+1}-\bt_nR_{n-1}=\frac{1-\al_n^2}{2},
\eea
which has the ``integrating factor''$R_n$. Therefore,
\bea
\mbox{} \bt_nR_nR_{n-1}&=&C_2+
\frac{1}{2}\sum_{j=0}^{n-1}\left(1-\al_j^2\right)R_j\nonumber\\
\mbox{} &=&C_2+\frac{1}{2}\sum_{j=0}^{n-1}
\left(1-\frac{C_1^2}{(4R_j^2-1)^2}\right)R_j,
\eea
where $C_2$ is another integration constant to be determined by
the initial condition
$$\bt_1=\frac{h_1}{h_0}=\frac{h_1}{\mu_0}=\frac{\mu_2}{\mu_0}-
\left(\frac{\mu_1}{\mu_0}\right)^2=\frac{4(\al+1)(\bt+1)}
{(\al+\bt+2)^2(\al+\bt+3)}.$$
After some computations,
$$C_2=\bt_1R_0R_1-\frac{1}{2}(1-\al^2_0)R_0=0.$$
Now the sum (2.17), may look complicated, however, with a 
partial fraction expansion, the sum can be taken and 
leads to
\bea
\mbox{} \bt_n = \frac{n(n+\al)(n+\bt)(n+\al+\bt)}
{(2n+\al+\bt)^2R_nR_{n-1}}.\nonumber
\eea
Therefore, after some simplifications we establish  
\bea
\mbox{} 
\b_n = \frac{4n(n+\al)(n+\bt)(n+\al+\bt)}
{(2n+\al+\bt)^2(2n+\al+\bt+1)(2n+\al+\bt-1)}.
\eea

\setcounter{equation}{0}
\setcounter{thm}{0} 

\section{Explicit Formulas}

We first determine ${\cal P}_n^{(\al, \b)}(\pm
1)$. Write (1.3) as,
\bea
\begin{gathered}
\frac{d}{d z}\,  {\cal P}_n^{(\al, \b)}(z) =
\frac{(n+r_n) {\cal P}_n^{(\al, \b)}(z) - \bt_n R_n
{\cal P}_{n-1}^{(\al,
\b)}(z)}{z-1} \\
\qquad \qquad \quad  + \frac{\bt_n R_n
{\cal P}_{n-1}^{(\al, \b)}(z) - r_n {\cal P}_n^{(\al, \b)}(z)}{z+1},
\end{gathered} 
\nonumber
\eea
and since $\frac{d}{d z}\,  {\cal P}_n^{(\al, \b)}(z)$ is regular at
$z=\pm 1,$ we arrive at
\bea
(n+r_n){\cal P}_n^{(\al, \b)}(1)-\bt_n {\cal R}_n P_{n-1}^{(\al,
\b)}(1)=0,\nonumber
\eea
\bea
\bt_n R_n {\cal P}_{n-1}^{(\al, \b)}(-1) - 
r_n {\cal P}_n^{(\al, \b)}(-1)=0.\nonumber
\eea
Thus we find 
\bea
 {\cal P}^{(\al, \b)}_n(1) =
{\cal P}^{(\al, \b)}_0(1) \prod_{j=1}^{n}\frac{\bt_j R_j} {r_j+j},
\eea
and
\bea
  {\cal P}^{(\al, \b)}_n(-1)=
{\cal P}^{(\al, \b)}_0(-1)\prod_{j=1}^{n}\frac{\bt_jR_j}{r_j}.
\nonumber
\eea
Substituting for $\b_n$ $r_n$ and $R_n$ from (2.9), (2.12), 
and (2.18), and applying (2.15) we prove that 
\bea
  {\cal P}^{(\al, \b)}_n(-1) = \frac{(-1)^n}{2^n} \prod_{j=1}^{n}
\frac{(j+\bt)(j+\al+\bt)} {\left[j + (\al+\bt/2)\right] 
\left[j +(\al+\bt-1/2)\right]}. \nonumber 
\eea
Using the facts $(\l)_n = \Gamma(\l+n)/\Gamma(\l)$, $(2\l)_{2n} 
= 4^2 (\l)_n(\l+1/2)_n$ we rewrite the above equation as 
\bea
 {\cal P}^{(\al, \b)}_n(-1) = 
\frac{(-1)^n\; 2^n\; (\b+1)_n}{(\al+\b+n+1)_n}.
\eea
Similarly
\bea
{\cal P}_n^{(\al, \b)}(1)  =
\frac{ 2^n\; (\al+1)_n}{(\al+\b+n+1)_n}.
\eea

We next evaluate  $h_n$, the squares of the $L^2$ norms. In general
(1.1) and (1.2) yield, \cite{Rai}
\bea
h_n = h_0 \b_1\b_2 \cdots \b_n.
\eea
The beta integral evaluation gives 
\bea
h_0 = 2^{\al+\b+1}\,\frac{\Gamma(\al +1)\Gamma(\b+1)}
{\Gamma(\al+\b+1)}.
\eea
Thus 
\bea
\int_{-1}^1 {\cal P}_m^{(\al, \b)}(x) {\cal P}_n^{(\al, \b)}(x)
(1-x)^\al(1+x)^\b dx = h_n \delta_{m,n},
\eea
with
\bea
h_n=  \frac{2^{\al+\b +n +1}\,\Gamma(\al + n +1)
\Gamma(\b+1) \, n!}{(\al+\b + n + 1)_{n}\, 
\Gamma(\al +\b + n +2n+2)}. 
\eea

We now prove that 
\bea
\frac{d}{dz}{\cal P}_n^{(\al, \b)}(z) = n {\cal P}_n^{(\al+1,
\b+1)}(z)
\eea
For $\al > -1$, $\b > -1$,  and $m <n-1$ integration by parts gives
\bea
\begin{gathered}
\int_{-1}^1 x^m \left(\frac{d}{dx}{\cal P}_n^{(\al, \b)}(x)\right)
(1-x)^{\al+ 1}(1+x)^{\b+1} \, dx \\
=-  \int_{-1}^1 {\cal P}_n^{(\al, \b)}(x)\;  f(x)
(1-x)^{\al}(1+x)^{\b }\, dx
\end{gathered}
\nonumber 
\eea
where $f(x) = x^{m-1}[m +
x(\b-\al) - x^2(\al +\b +m+2)]$. Since $f$ has degree at most
$n-1$, the above integral must vanish and we conclude that 
$\frac{d}{dx}{\cal P}_n^{(\al, \b)}(x)$ is orthogonal to all
polynomials of degree less than $n-1$ with respect to $w^{(\al+1,
\b+1)}(x)$. The uniqueness of the orthogonal polynomials and the fact
that ${\cal P}_n^{(\al, \b)}(x), n \ge 0$ are monic, establish
(3.8). Clearly (3.8) and (3.2) give 
\bea
\begin{gathered}
\left. \frac{d^k}{dx^k} {\cal P}_n^{(\al, \b)}(x) \right|_{x=-1} =
\frac{n!}{(n-k)} {\cal P}_{n-k}^{(k+\al, k+\b)}(-1) \\ 
\mbox{} \qquad \qquad  \qquad =
\frac{(-2)^{n-k}(\b+k)_{n-k}}{(\al+\b + n+k+1)_{n-k}}
\end{gathered}
\eea
The Taylor series about $x=-1$ now gives the representation
\bea
\begin{gathered}
{\cal P}_n^{(\al, \b)}(x) =
\frac{(-2)^{n}(\b+1)_{n}}{(\al + \b +1)_{n}}  \\
\mbox{} \qquad \qquad \times {}_2F_1(-n,
n+\al+\b +1; \b+1; (1+x)/2),
\end{gathered}
\eea
which we recognized to be the monic Jacobi polynomials.
Similarly (3.3) and (3.8) give the alternate representation 
\bea
\begin{gathered}
{\cal P}_n^{(\al, \b)}(x) =
\frac{(2)^{n}(\al +1)_{n}}{(\al + \b +1)_{n}}  \\
\mbox{} \qquad \qquad \times {}_2F_1(-n,
n+\al+\b +1; \al + 1; (1- x)/2),
\end{gathered}
\eea

\end{document}